\documentclass{mrlart3}
  \def\volno{14}
  \def\yrno{2007}
  \def\issueno{3}
  \def\lpageno{511} 
\usepackage{amsmath}
\usepackage{amssymb}
\usepackage{amsthm}
\usepackage{oldgerm}

\newcommand{\Rm}{\mathbb{R}}

\newcommand{\Nm}{\ensuremath{\mathbb{N}}}

\newcommand{\mA}{\ensuremath{\mathcal{A}}}

\newcommand{\mI}{\ensuremath{\mathcal{I}}}

\newcommand{\Tm}{\ensuremath{\mathbb{T}}}

\newtheorem{lem}{Lemma}
\newtheorem{thm}{Theorem}

\newtheorem{prop}[lem]{Proposition}
\newtheorem{defn}[lem]{Definition}

\def\lto{\longrightarrow}
\def\lmto{\longmapsto}

\def\leq{\leqslant}
\def\geq{\geqslant}

\address{Patrick Bernard\\
CEREMADE, UMR CNRS 7534\\
Universit\'e de Paris Dauphine\\
Pl. du Mar\'echal de Lattre de Tassigny\\
75775 Paris Cedex 16\\
France}

  \email{patrick.bernard@ceremade.dauphine.fr}
  \urladdr{http://www.ceremade.dauphine.fr/ pbernard/}

\title{Smooth critical sub-solutions of the Hamilton-Jacobi equation}
\author{Patrick  Bernard}
\date{  2006 }
\begin{document}

\maketitle
\vspace{1cm}
\begin{abstract}
We establish the existence of smooth critical sub-solutions
of the Hamilton-Jacobi equation on compact manifolds for smooth
convex Hamiltonians, that is in the context of weak KAM theory,
under the assumption that the Aubry set is the union of finitely
many hyperbolic  periodic orbits or fixed points.
\end{abstract}

\section*{}
Let $M$ be a compact manifold without boundary.
A function $H(x,p):T^*M\lto \Rm$ is called a Tonelli
Hamiltonian if it is $C^2$ and if, 
 for each $x\in  M$
 the function $p\lmto H(x,p)$
 is convex with positive definite Hessian and superlinear
 on the fibre $T_x^*M$.
 A Tonelli Hamiltonian generates a complete Hamiltonian flow $\psi_t$.
It is known that there exists one and only one real constant
$\alpha(H)$  such that the equation 
\begin{equation}\tag{HJ}
H(x,du_x)=\alpha(H)
\end{equation}
has a solution in the viscosity sense.
This equation, then, may have several solutions.
These solutions are  Lipschitz  but, in general,
none of them is  $C^1$.
It is a natural question whether there exist sub-solutions
of (HJ) which are more regular.
In this direction, Fathi and Siconolfi proved the
existence of  a $C^1$ sub-solution, see 
\cite{FS:04}.
More recently, I proved the existence of $C^{1,1}$ 
sub-solutions  see \cite{inv}.
Examples show that, in general, 
$C^2$ sub-solutions do not exist.
In the present paper, we establish, under the additional
assumption that the Aubry set is a finite
union of hyperbolic orbits and hyperbolic fixed points,
the existence of smooth sub-solutions.
This answers a question asked several times to the author
and to Albert Fathi during various conferences
on weak KAM theory.
This also provides a nice class of examples, and is 
a useful technical step for deeper studies of this class
of examples, see \cite{AIPSM} and other papers.

\section{introduction}
Let $M$ be a compact Manifold and $H$ a Tonelli
Hamiltonian on $T^*M$.
We denote by $\pi:T^*M\lto M$ the canonical projection.
We say that the function 
$u: M\lto \Rm$
is a sub-solution of the equation
\begin{equation}\tag{HJc}
H(x,du_x)=c
\end{equation}
if it is Lipschitz continuous and  satisfies the inequality
$
H(x,du_x)\leq c
$
at all its points of differentiability.
It is equivalent to say that $u$ is a sub-solution
in the viscosity sense
(see \cite{Fa:un,Ba:94,CaSi} and many other texts for a definition), or that $u$ is Lipschitz and that the inequality above holds 
at Lebesgue almost every point.
It is known, see \cite{CIPP} for example, that 
the set of real numbers $c$ such that a sub-solution exists
is of the form $[\alpha(H),\infty)$, for a given real number
$\alpha(H)$.
We take this as a definition of $\alpha(H)$, and call this number
the Ma\~n\'e critical value.
If $c>\alpha(H)$, then the equation $(HJc)$ has smooth sub-solutions,
see \cite{CIPP}.
It is known that $\alpha(H)$ is the only value of $c$ such that 
the equation $(HJc)$ has a  solution in the viscosity sense,
see \cite{Fa:un}, or \cite{FS:04}.

A sub-solution of (HJ) is called strict at $x$ if there exists a 
neighborhood $U$ of $x$ and a positive number $\delta$
such that the inequality
$$
H(x,du_x)\leq \alpha(H)-\delta
$$
holds at all point of differentiability $x$
of the function $u$
in $U$.

One of the main reasons why we pay a special attention to 
sub-solutions is that they are calibrators for the
problem of minimization of the Lagrangian action.
In order to be more precise, we associate to the Hamiltonian 
$H$ the Lagrangian $L: TM\lto \Rm$
defined by
$$
L(x,v)=
\max_{p\in T_q^*M}\big(
p(v)-H(x,p)\big).
$$
The function $u:M\lto \Rm$ 
 is a sub-solution of (HJc)
if and only if the inequality
\begin{equation}\label{calibration}
\int_S^T
c+ L(\gamma(t),\dot\gamma) dt
\geq
u(\gamma(T))-u(\gamma(S))
\end{equation}
holds for all absolutely continuous curves 
$\gamma:[S,T]\lto M$ (see \cite{Fa:un}).
Of course, having equality in this inequality implies
that the curve $\gamma$ is minimizing the action with fixed
endpoints.
Following Fathi, we say that the
absolutely continuous  curve $\gamma:I\lto M$,
where $I\subset \Rm$ is an interval,
is calibrated by $u$ (for $(H,c)$)
if there is equality in (\ref{calibration})
for all $[S,T]\subset I$.
This implies that $\gamma$ is $C^2$
and that there exists a Hamiltonian trajectory
$(\gamma(t),p(t))$ above $\gamma$.

It is known, see \cite{Fa:un},
that there exist curves $\gamma\in C^2(\Rm, M)$
calibrated by $u$ for $(H,c)$ on the whole real line
if and only if $c=\alpha(H)$.
Given a sub-solution $u$ of (HJ),
we define the nonempty compact invariant set
$\tilde \mI(H,u)$
as the union of the images of all the Hamiltonian trajectories
of $H$ whose projection on $M$ are calibrated by $u$ on the whole 
real line.
We then define the Aubry set 
$$
\tilde \mA(H):=
\bigcap_u \tilde \mI(H,u),$$
where the intersection is taken on all sub-solutions of (HJ).
It is known that this set is not empty.
We denote by $\mI(H,u)$ and $\mA(H)$
the projections of $\tilde \mI(H,u)$ and 
$\tilde \mA(H)$ on $M$. these are compact and non-empty 
subsets of $M$.
We are now ready to state our main result:

\begin{thm}\label{main}
Let $H$ be a $C^k$ Tonelli Hamiltonian, $2\leq k\leq \infty$.
Assume that the Aubry set $\tilde \mA(H)$ is a finite union
of hyperbolic periodic orbits and hyperbolic fixed points.
Then there exists a $C^k$ sub-solution of (HJ)
which is strict outside of the projected 
Aubry set $\mA(H)$.
\end{thm}

\section{Examples}\label{examples}
\subsection{Mechanical Hamiltonian system}
Let us consider the case
$$
H(x,p)=\frac{1}{2}\|p\|_x^2+V(x)
$$
where $\|.\|_x$ is a Riemannian metric on $M$
and $V$ is a smooth function on $M$.
Then it is easy to see that 
$\alpha(H)=\max V$,
and that there exists a smooth sub-solution to (HJ):
any constant function is such a sub-solution!

\subsection{Non-existence of a $C^2$ sub-solution}
Let us now specialize to $M=\Tm$,
and consider the Hamiltonian 
$$
H_P(x,p)=\frac{1}{2}(p+P)^2-\sin^2 (\pi x)
$$
depending on the real parameter $P$.
For $P=0$ this is a Mechanical system as discussed above,
and the constants are sub-solutions of (HJ).
Let $X(x):\Tm\lto \Rm$ be the function
such that $X(x)=\sin(\pi x)$ for $x\in [0,1]$.
For $P=0$, there is  one and only one 
viscosity solution, and this solution 
is not differentiable.
This solution is the primitive of the function 
$s(x)X(x)$, where $s(x)=1$ on $[0,1/2[$
and $s(x)=-1$ on $[1/2,1[$.

Let us set
$$
a=\frac{2}{\pi}
=\int_{\Tm} X(x) dx.
$$
For each smooth function $f:\Tm\lto M$
let us denote by $P_f$ the average of $f$, and by $F$
a primitive of $f-P_f$.
It is not hard to see that, if we assume that 
$|f(x)|\leq X(x)$, then 
the function $F$ is a sub-solution 
of $(HJ)$ for the value $P_f$ of the parameter.
We conclude that, for $P\in [-a,a]$,
we have $\alpha(H_P)=0$.
In addition, for each $P\in ]-a,a[$,
the equation (HJ) has smooth sub-solutions.
For these values of $P$, the Aubry set is the fixed point $(0,0)$.

However, for $P=a$, there is one and only one sub-solution
of (HJ), which turns out to be a solution.
It is given by the primitive of the function 
$X-a$.
This function is $C^{1,1}$  but not $C^2$.
Note that the Aubry set, then, is not reduced to the 
hyperbolic fixed point $(0,-a)$ but is the whole graph of $X-a$.

\section{Some preliminaries}

\subsection{On Weak KAM Theory}
We recall some important facts of 
Weak KAM theory and  Aubry-Mather theory, see
 \cite{Ma:91, Ma:97,CDI, Fa:un} for the proofs.
Let us begin with the definition of viscosity solutions.

\begin{defn}
The function $u$ is a viscosity solution of
(HJc) if $u$ is a sub-solution and if,
for each $x\in M$, 
there exists a curve
$\gamma_x\in C^2((-\infty,0],M)$
which satisfies $\gamma(0)=x$ and which is calibrated by $u$
for $(H,c)$.
\end{defn}

This definition, due to Fathi, is the one we will use.
Fathi proves in \cite{Fa:un} that it is equivalent
to the standard definition of viscosity solutions.
Recall that there exists a viscosity solution
of $(HJc)$ if and only if
$c=\alpha(H)$.

An important observation of Weak KAM theory is that,
if $\gamma(t):[S,T]\lto M$ is calibrated by a sub-solution $u$
for (H,c), then the function $u$ is differentiable 
on $\gamma(]S,T[)$.
In addition, if $u$ is differentiable at 
$\gamma(t), t\in[S,T]$, then 
$du_{\gamma(t)}=p(t)$, where $p(t)$ is such that 
$(\gamma(t),p(t))$ is the orbit of the Hamiltonian flow of $H$
which projects on $\gamma(t)$.

 The basic result of Aubry-Mather theory
 and weak KAM theory is  that there exists a Lipschitz section 
 $X(x):\mA(H)\lto T^*M$ over $\mA(H)$ such that 
 $\tilde \mA(H)=X(\mA(H))$.
 In other words, if $(\gamma(t),p(t)):\Rm\lto T^*M$ is a Hamiltonian
 trajectory, then it is equivalent to say that the curve $\gamma$
 is calibrated by all sub-solutions of (HJ)
 and to say that $\gamma(0)\in \mA(H)$ and 
 $p(0)=X(\gamma(0))$.
 All the sub-solutions of (HJ) are differentiable on 
 $\mA(H)$, and satisfy
 $du_x=X(x)$ at each point $x\in \mA(H)$.

\subsection{On hyperbolic periodic orbits and fixed points}
We say that the fixed point $P$ is a Hyperbolic fixed
point if the linearized equation at $P$ does not
have any eigenvalue of zero real value.
In this case, the tangent space $T_P(T^*M)$
splits as the sum of two Lagrangian subspaces
$L+$ and $L-$, the stable and unstable spaces,
which are invariant by the linearized map, and such that
the linearized flow on $L^+$ has only eigenvalues
of positive  real part,
while the linearized flow on  $L^-$
has eigenvalues of negative real part.

Let now $P(t):\Rm\lto T^*M$ be a $T$- periodic orbit.
We denote by $\tilde A$ the image of this orbit, 
which is a simple closed
curve in $T^*M$.
For each point $P\in \tilde A$, we consider the linear 
endomorphism 
$d\psi_T(P)$ of  $T_P(T^*M)$.
This endomorphism has $1$ as an eigenvalue with multiplicity 
at least two.
We say that the orbit $\tilde A$ is hyperbolic if,
for each $P\in \tilde A$, the multiplicity of $1$
as an eigenvalue of $d\psi_T(P)$ is exactly $2$,
and if the other eigenvalues have modulus different from one.
Then, there exist,
at each point $p\in \tilde A$
two Lagrangian subspace $W^+(P)$ and $W^-(P)$
of $T_P(T^*M)$, which are invariant under $d\psi_T(P)$
and such that the restriction of this linear map 
to $W^+(P)$ has all but one of its  eigenvalues of modulus
greater than  one,
and the restriction to $W^-(P)$ has all but one of its eigenvalues
of modulus smaller than one.
These subspaces depend continuously on $P$.
In addition, the spaces $W^+(P)$ and $W^-(P)$
are transversal in the kernel of $dH_P$,
which means that 
$W^+(P) +W^-(P)=\ker(dH_P)$.
 Finally, we have
$W^+(P)\cap W^-(P)=\Rm Y(P)$, where $Y$ is the Hamiltonian
vectorfield.

If a hyperbolic fixed point or a hyperbolic periodic orbit
is minimizing, in the sense that the orbit
$\pi \circ P(t): \Rm\lto M$
is minimizing the action with fixed endpoints on each 
compact interval, then the stable and unstable subspaces $W^{\pm}$
are transversal to the vertical at each point.
This is because they have to coincide with
the so-called Green spaces, see \cite{CoIt}, Proposition B.
We then have:

\begin{thm}
If $\tilde A\subset \tilde \mA$ is either a 
hyperbolic fixed point or
a hyperbolic periodic orbit, it has  a $C^{k-1}$
local unstable  manifold $W^+$, which is locally  the graph 
of a $C^{k-1}$ closed  one-form. 
There exists a neighborhood
$\tilde V$ of $\tilde A$ such that
each point $P\in T^*M$ whose backward orbit is contained in 
$\tilde V$ satisfies $P\in W^+$.
\end{thm}

\section{Outline of the proof}
Of central interest will be the notion 
of \textit{sharp potential}.

\begin{defn}
A sharp potential is a smooth function 
$V: M\lto \Rm$
which is null on $\mA(H)$ and positive outside of this set.
\end{defn}

Given a function (and in particular a sharp potential)
$V: M\lto \Rm$,
we denote by $H+V$ the Hamiltonian
$$
(H+V)(x,p)=H(x,p)+V(x).
$$
We will prove in section \ref{potentials}
the following general result:

\begin{thm}\label{thmpot}
Assume that $H$ is a $C^2$ Tonelli Hamiltonian, 
then there exists a sharp potential $V_0$ such that,
for all sharp potential $V\leq V_0$, we have
$\alpha(H+V)=\alpha(H)$
and $\tilde \mA(H+V)=\tilde \mA(H)$.
\end{thm}

Assume now that the Aubry set $\mA(H)$
is the union of finitely many hyperbolic periodic orbits
or fixed points.
Then the Aubry set $\mA(H+V)$
is the union of the same family of points and closed
curves, which are also Hamiltonian orbits of
$H+V$.
In addition, the potential $V$ in Theorem \ref{thmpot}
 can be chosen very flat on 
$\mA(H)$, in such a way that the linearized Hamiltonian
flow along these orbits is the same for 
$H$ and $H+V$.
As a consequence, the orbits remain hyperbolic as orbits 
of $H+V$.
We will prove in Section \ref{local}:

\begin{thm}\label{localthm}
Let $H$ be a $C^k$ Tonelli Hamiltonian, $2\leq k\leq \infty$.
Assume that the Aubry set $\tilde \mA(H)$ is a finite union
of hyperbolic periodic orbits or fixed points.
Then there exists a solution of (HJ) which is 
$C^k$ in a neighborhood of $\mA(H)$.
\end{thm}

Applying Theorem \ref{localthm} to the Hamiltonian 
$H+V$, we obtain a solution of the Hamilton-Jacobi
equation 
$$
H(x,du_x )+V(x)=\alpha(H+V)
=\alpha(H)
$$
which is smooth in a neighborhood of $\mA(H+V)=\mA(H)$.
This function is then a sub-solution of (HJ)
which is strict outside of the Aubry set
and smooth in a neighborhood of the Aubry set.
It is easy to see  that this function can be smoothed out
to a $C^k$ sub-solution of (HJ).
We have reduced the proof of Theorem \ref{main}
to the proof of Theorems \ref{thmpot} and \ref{localthm}.
Note that smoothing in this kind of context has been 
done in \cite{CIPP}, \cite{FM}, Theorem 8.1,
and \cite{FS:04}, Theorem 9.2.

\section{Sharp potentials and strict  sub-solutions}\label{potentials}
This section is general. 
We work with an arbitrary Tonelli Hamiltonian
and make no assumption on it's Aubry set.
We prove Theorem \ref{thmpot}.

\begin{prop}\label{equiv}
The following statements are equivalent:
\begin{enumerate}
\item
There exists a Lipschitz sub-solution which is strict at
each point $x$ of the complement of $\mA(H)$.
\item
There exists a sharp potential  $V_0:\Tm\times M\lto \Rm$,
such that,
for each sharp potential  $V\leq V_0$, we have 
$\alpha(H+V)=\alpha(H)$.
\item
There exists a sharp potential  $V_0:\Tm\times M\lto \Rm$,
such that,
for each sharp potential  $V\leq V_0$, we have 
$\alpha(H+V)=\alpha(H)$ and $\tilde \mA(H+V)=\tilde \mA(H)$.
\end{enumerate}
\end{prop}

\begin{proof}
It is enough to prove that 
$2\Rightarrow 1\Rightarrow 3$. 
Let us first assume 2.
Let $u$ be a 
Lipschitz sub-solution of the Hamilton-Jacobi equation
\begin{equation}\tag{HVJ}
H(x,du_x) +V(x)=\alpha(H+V)
=\alpha(H).
\end{equation}
Such a solution exists by the definition of $\alpha(H+V)$.
Then it is clear that $u$ is also a Lipschitz sub-solution
of (HJ), which is strict outside of $\mA(H)$.

Let us now assume 1.
We have a Lipschitz sub-solution 
$u$ which is strict outside of $\mA(H)$.
Then, there exists a sharp potential $V_0$ such that
$2V_0(x)+H(x,du_x)\leq \alpha(H)$ at all point of 
differentiability of $u$.
As a consequence, for each
sharp potential $V\leq V_0$, the function $u$ is a sub-solution of the equation 
$$
H(x,du_x) +V(x)
=\alpha(H)
$$
which is strict outside of $\mA(H)$.
This implies that $\alpha(H+V)=\alpha(H)$.
In addition, the sub-solution $u$ of (HVJ)
being strict outside of $\mA(H)$, there is no point of
$\mA(H+V)$ outside of $\mA(H)$, or in other words
$\mA(H+V)\subset \mA(H)$.
We now claim that  $\tilde \mA(H+V)\supset \tilde \mA(H)$,
which terminates the proof of the equality between the Aubry sets.
Let us consider an orbit  $(q(t),p(t)):\Rm\lto T^*M$
of the Hamiltonian flow of $H$ which is contained in $\tilde \mA(H)$.
This curve is also an orbit of
the Hamiltonian flow of  $H+V$, because $V$ is flat on 
$\mA(H)$.
In order to prove that this orbit is contained in 
$\tilde \mA(H+V)$,
it is enough to prove that the curve $q(t)$
is calibrated by all  sub-solutions of
(HVJ).
 This is true because the curve $q(t)$ is calibrated by
 all sub-solutions of (HJ), and because each sub-solution
of (HVJ) is a sub-solution of 
 of (HJ).
\end{proof}

This is relevant in view of:

\begin{prop}\label{true}
The equivalent properties of Proposition
\ref{equiv}  hold true for all Tonelli
Hamiltonians.
\end{prop}
\begin{proof}
The item 1 has been established by Fathi and Siconolfi
in \cite{FS:04}.
It was one of the key steps in the construction of a $C^1$
sub-solution.
Another proof
 consists of establishing directly item 2.
Let us explain this other proof (which is close to 
the one of Massart, see \cite{Ma:un}).
We shall first state without proof a 
Lemma which was introduced 
in an easier  situation in \cite{BeBu:first},
and has been proved by 
 Daniel Massart \cite{Ma:un} in the present setting:

\begin{lem}
Let $U$ be an open set whose closure is disjoint from
$\mA(H)$.
Then there exists a smooth function 
$W: M\lto \Rm$
which is positive on $U$ and null outside of $U$,
and such that $\alpha(H+aW)=\alpha(H)$ for all real numbers
$a\in [0,1]$.
\end{lem}

Let us now cover the complement of $\mA(H)$ by a countable family
$U_n$ of open sets whose closure do not intersect $\mA(H)$.
Let $W_1$ be the function given by the Lemma applied on $U_1$.
Let $W_2$ be the function given by the Lemma applied to the 
Hamiltonian $H+W_1$ on $U_2$.
We define by recurrence  the function $W_n$
as the function given by the Lemma applied to the Hamiltonian
$H+W_1+\cdots +W_{n-1}$ on $U_n$.
Then for each sequence $a_n$ of real numbers in $]0,1]$,
and each $n\in \Nm$,
we have 
$$
\alpha(H+a_1W_1+\cdots+a_nW_n)=\alpha(H).
$$
We can assume that the sequence $a_n$ is decreasing so fast
that the sequence  
$a_1W_1+\cdots +a_nW_n$ is converging uniformly
to a limit $W$, which is a continuous function
positive outside of $\mA(H)$.
Let then $V$ be a sharp potential such that 
$0\leq V\leq W$. We claim that $\alpha(H+V)=\alpha(H)$.
Indeed, for each $\epsilon>0$, there exists 
$n\in \Nm$ such that 
$$\epsilon +a_1W_1+\cdots+a_nW_n\geq W\geq V$$
but then we have, setting $S_n=a_1W_1+\cdots+a_nW_n$
$$
\epsilon+ \alpha(H)
=
\epsilon+
\alpha(H+S_n)=
\alpha(\epsilon +H+S_n)
\geq
\alpha(H+V)\geq
\alpha(H).
$$
Since this holds for all $\epsilon>0$, we obtain
$\alpha(H+V)=\alpha(H)$ as desired.
\end{proof}

All the content of the present section
remains true in  the case of time-periodic Hamiltonian 
systems.
The proof we present of Proposition \ref{true}
can be transposed easily to this more general setting, see
\cite{Ma:un},
while the original proof of Fathi and Siconolfi can not.

\section{Existence of a locally smooth solution}\label{local}
We prove Theorem \ref{localthm}.
In order to do so, we consider a $C^k$ Tonelli Hamiltonian 
$H$, and assume that the Aubry set $\tilde \mA(H)$
is a finite union of 
hyperbolic periodic orbits or fixed points.
Let us first recall, without proof, a useful result of Fathi,
\cite{Fa:un}.

\begin{prop}\label{I=A}
There exists a viscosity solution $u$ of (HJ)
such that $\tilde \mI(H,u)=\tilde \mA(H)$.
\end{prop}

The definition of the set $\tilde \mI(H,u)$
is given in the introduction.
Let us fix from now on one of the viscosity solutions
$u$ given by Proposition \ref{I=A},
and  prove that this viscosity solution satisfies 
the conclusions of Theorem \ref{localthm}.
Let $\tilde A$ be one of the  connected components of 
$\tilde \mA(H)$, and $A$ be its projection on $M$.
Note that both $A$ and $\tilde A $ are either points
or simple closed curves.
It is enough to prove:

\begin{prop}
The function $u$ is $C^k$ in a neighborhood of $A$.
\end{prop}

\begin{proof}
Let us call
$\Gamma_u$ the set of points 
$(q,p)\in T^*M$
such that the curve 
$q(t)=\pi\circ \psi_t(q,p)$ (where $\pi$ is the
projection on the base and $\psi_t$ is the Hamiltonian flow)
is calibrated by $u$ on $(-\infty,0]$.
It is known that the function $u$ is differentiable
at a point $x$ if and only if the set 
$\Gamma_u\cap T_x^*M$ contains only one point $p$, and then
$d_xu=p$.
We claim that the
germ of $\Gamma_u$ along $A$ is equal to the
germ of the local unstable manifold $W^+$.
It means that there exists a neighborhood
$ V$ of $\tilde A$ such that 
$\Gamma_u\cap  V=W^+\cap  V$.
Since $W^+$ is the graph of a $C^{k-1}$ one-form,
this implies the desired result.
The claim follows from:

\begin{lem}\label{orbites}
Let $ V$ be a compact neighborhood of $\tilde A$ satisfying 
$ V\cap \tilde \mA(H)=\tilde A$.
There exists a neighborhood $U$ of $A$ such that
 all the curves 
$\gamma \in C^2((-\infty,0],M)$
which are calibrated by $u$
and satisfy $\gamma(0)\in U$
satisfy $(\gamma(t),\xi(t))\in  V$ for all $t\in (-\infty,0] $,
where $(\gamma(t),\xi(t))$ is the Hamiltonian trajectory lifting
$\gamma$.
\end{lem}

\begin{proof}
Assume that the conclusion does not hold.
Then, there exists a sequence $(\gamma_n(t),\xi_n(t))$
of trajectories in $C^1((-\infty,0],T^*M)$
such that the curves $\gamma_n$ are calibrated on $(-\infty,0]$,
and such that  
$$\lim _{n\lto \infty}d(\gamma_n(0), A)=0,$$
and  a sequence $-T_n\leq 0$ of times  such that 
$\gamma_n(-T_n)\in \partial V$ (the boundary of $V$).
Let $(q_n,p_n):(-\infty,T_n]\lto T^*M$
be the trajectory $(q_n(t),p_n(t))=(\gamma_n(t-T_n),\xi_n(t-T_n))$.
We can assume by taking a subsequence that the sequence
$(q_n(t),p_n(t))$ converges uniformly on compact sets to a limit 
trajectory
$(q(t),p(t)):I\lto M$, where $I$ is either an interval
of the form $(-\infty,T]$ or is $\Rm$.
It is easy to see that the curve $q(t)$ is calibrated by $u$,
and that $(q(0),p(0))\in \partial V$, hence 
$(q(0),p(0))\notin \tilde \mA(H)$.

If $I=(-\infty,T]$, then we have $q(T)\in A$
and, since $u$ is differentiable on $A$, we have
$p(T)=du_{q(T)}=X(q(T))$, where
 $X(q(t))=\tilde \mA(H)\cap T_{q(t)}^*M$.
But this implies that the curve $(q(t),p(t))$
is contained in $\tilde \mA(H)$, which is 
a contradiction.

If $I=\Rm$, then $(q(0),p(0))\in \tilde \mI(H,u)$,
which is in contradiction with the fact that
$\tilde \mI(H,u)= \tilde \mA(H)$.
\end{proof}
This ends the proof of Theorem \ref{localthm}. \end{proof}

\section*{Acknowledgements}
This paper was born from conversations which occurred
during the Workshop on Conservative Dynamics in Rio, august 2005.
I wish to thank Leonardo Marcarini for organizing this
meeting, and also for participating in these conversations and
encouraging me to write the present paper.
I also thank Albert Fathi for listening patiently
to my successive attempts of proving the existence of smooth 
sub-solutions.
Gonzalo Contreras and Renato Iturriaga were also 
interested and motivating interlocutors, they proposed some
simplifications in the proofs from which the present text takes 
advantage.

\small
\bibliographystyle{amsplain}

\end{document}